\documentclass[11pt, twoside]{article}
\usepackage{amssymb, amsmath}
\usepackage[matrix,arrow]{xy}
\usepackage{amscd}

\begin{document}
\newtheorem{Mthm}{Main Theorem.}
\newtheorem{Thm}{Theorem}[section]
\newtheorem{Prop}[Thm]{Proposition}
\newtheorem{Lem}[Thm]{Lemma}
\newtheorem{Cor}[Thm]{Corollary}
\newtheorem{Def}[Thm]{Definition}
\newtheorem{Guess}[Thm]{Conjecture}
\newtheorem{Ex}[Thm]{Example}
\newtheorem{Rmk}{Remark.}
\newtheorem{Not}{Notation.}

\renewcommand{\theThm} {\thesection.\arabic{Thm}}
\renewcommand{\theProp}{\thesection.\arabic{Prop}}
\renewcommand{\theLem}{\thesection.\arabic{Lem}}
\renewcommand{\theCor}{\thesection.\arabic{Cor}}
\renewcommand{\theDef}{\thesection.\arabic{Def}}
\renewcommand{\theGuess}{\thesection.\arabic{Guess}}
\renewcommand{\theEx}{\thesection.\arabic{Ex}}
\renewcommand{\theRmk}{}
\renewcommand{\theMthm}{}
\renewcommand{\theNot}{}
\renewcommand{\thefootnote}{\fnsymbol{footnote}}

\newcommand{\spec}{\textnormal{Spec}\hspace{1mm}}
\newcommand{\proj}{\textnormal{Proj}\hspace{1mm}}
\newcommand{\Diff}{\textnormal{Diff}\hspace{1mm}}
\newcommand{\pf}{{\bfseries\itshape Proof. }}
\newcommand{\mult}{\textnormal{mult}}
\newcommand{\com}{\hspace{-2mm}\textnormal{\textbf{.}}\hspace{2mm}}
\newcommand{\bir}{-\hspace{-1mm}\rightarrow}
\newcommand{\qed}{\hfill Q.E.D.\newline}
\phantom{this is space} \vspace{46mm} \Large
\textbf{\hspace{-6mm}Log Canonical Thresholds and Generalized
Eckardt Points}
\vspace{3mm}\\
\normalsize
\textbf{Ivan Cheltsov and  Jihun Park}
\vspace{27mm}\\
\small {\textbf{Abstract.} Let $X$ be a smooth hypersurface of
degree $n\geq 3$ in $\mathbb{P}^n$. We prove that the
log canonical threshold of $H\in|-K_X|$ is  at least
$\frac{n-1}{n}$. Under the assumption of the Log minimal model
program, we also prove that a hyperplane section $H$ of
$X$ is a cone in $\mathbb{P}^{n-1}$ over a smooth hypersurface of
degree $n$ in $\mathbb{P}^{n-2}$ if and only if
the log canonical threshold of $H$ is $\frac{n-1}{n}$.\\
Bibliography :  20 titles.} \normalsize

\footnotetext{This work was partially supported by NSF Grant
DMS-0100991.}
\thispagestyle{empty}

\renewcommand{\thesection}{\large{\arabic{section}.}}
\section{\hspace{-3mm}\large{Introduction}}
\renewcommand{\thesection}{\arabic{section}}
All varieties are assumed to be defined over $\mathbb{C}$, unless
otherwise stated. Main definitions and notations appear in
\cite{KMM} and \cite {Sho92}.

To measure how far from log canonicity log pairs are,
V.~Shokurov introduced log canonical thresholds in \cite{Sho92}.
It is known that they have many amazing properties.
\begin{Def}\com
Let $(X,B)$ be a log canonical pair and let  $Z$ be a closed
subvariety of $X$. Suppose that $D$ is a $\mathbb{Q}$-Cartier
divisor on $X$. The log canonical threshold of $D$ along $Z$ with
respect to $K_X+B$ is the number:
\[lct_Z(D;X,B)=\sup\{c:K_X+B+cD \textnormal{ is log canonical along }Z\}.\]
\end{Def}
It is easy to check $0\leq lct_Z(D;X,B)\leq 1$. If $B=0$, we use
$lct_Z(D;X)$ instead of $lct_Z(D;X,0)$. For the case $Z=X$ we use
the notation $lct(D;X,B)$ instead of $lct_X(D;X,B)$.

We may find log canonical thresholds in several other branches of
mathematics in various disguises. For
example, we consider the log canonical threshold
$lct_0(D;\mathbb{C}^n)$ of $D=(f=0)$ at the origin with
respect to $(\mathbb{C}^n,0)$, where $f$ is a nonconstant
holomorphic function near the origin. Then we can see (\cite{Ko97}) that this
number is the same as the following number:
\[\sup\{c:|f|^{-c} \textnormal{ is locally }L^2 \textnormal{near the
origin} \}.\]

Bernstein-Sato polynomials provide another example. Bernstein-Sato
 polynomials appear in differential operator theory---in
 particular, $\mathcal{D}$-module theory (\cite{Bj79}). Let us briefly explain
 what Bernstein-Sato polynomials are. It is known that for any
 convergent power series $f\in \mathbb{C}\{z_1,\cdots,z_n\}$,
 there is a nonzero polynomial $b(s)\in\mathbb{C}[s]$ and a linear
 differential operator
 $Q=\sum_{I,j}f_{I,j}s^j\frac{\partial^I}{\partial z_I}$ such that
 \[b(s)f^s=Qf^{s+1},\]
 where each $f_{I,j}$ is a convergent power series.
 For a fixed $f$, such $b(s)$'s form an ideal of $\mathbb{C}[s]$.
 The monic generator of this ideal is called the Bernstein-Sato
 polynomial of $f$. We can see that $lct_0(D;\mathbb{C}^n)$
 is the absolute value of the largest root of the Bernstein-Sato
 polynomial of $f$ (\cite{Ko97}).

 Recently, M.~Musta\c{t}\v{a} investigated log canonical thresholds
via jet schemes (\cite{Mu00}). He obtained
\begin{Thm}\com
Let $X$ be a smooth variety and $D$ an integral effective divisor
on $X$. Then the log canonical threshold of $D$ with respect to $K_X$
is given by
\[lct(D;X)=\dim X-\sup_{m\geq 0}\frac{\dim D_m}{m+1},\]
where $D_m$ is the $m$-th jet scheme of $D$.
\end{Thm}
The $m$-th jet scheme $X_m$ of a variety $X$ is a scheme whose
closed points over $x\in X$ are morphisms
$\mathcal{O}_{X,x}\longrightarrow k[t]/(t^{m+1})$. If $X$ is a
smooth variety, then $X_m$ is an affine bundle over $X$ of
dimension $(m+1)\dim X$.
\vspace{5mm}

To understand a given variety, it is important to investigate
linear systems related to the canonical divisor. One  such
investigation is to find ``extreme" elements in the linear
systems. We have two kinds of extreme elements in  the linear
systems. One is a ``good" element, and the other is a ``bad"
element. We need to explain what ``good" elements are and what
``bad" elements are. It is natural that singularities should
distinguish between the ``good" and the ``bad". Since we always
consider linear systems related to canonical divisors (or log
canonical divisors), these concepts should involve
canonical divisors.

For a ``good" element, M.~Reid considered a general elephant. Following him,
V.~Shokurov introduced more general concepts.
\begin{Def}\com
Let $X$ be a normal variety and let $D=S+B$ be a subboundary on
$X$ such that $S$ and $B$ have no common components, $S$ is a
reduced divisor, and $\lfloor B\rfloor\leq 0$. We say that
$K_X+D$ is n-complementary if there is a divisor $D^+$ on $X$
satisfying the following conditions:
\begin{enumerate}
\item
$nD^+$ is integral and $n(K_X+D^+)$ is linearly trivial.
\item
$nD^+\geq nS+\lfloor (n+1)B\rfloor$.
\item
$K_X+D^+$ is log canonical.
\end{enumerate}
The divisor $K_X+D^+$ is called an $n$-complement of $K_X+D$.
\end{Def}

We now introduce a counterpart of these ``good" elements.
it was originally introduced by S.~Keel and J.~McKernan
(\cite{KeMc99}). Strictly speaking, it is a counterpart of 
M.~Reid's general elephant.
\begin{Def}\com
Let $X$ be a normal variety. Let $B$ be an effective
$\mathbb{Q}$-divisor on $X$. A special tiger for $K_X+B$ is an
effective $\mathbb{Q}$-divisor $D$ such that $K_X+B+D$ is numerically
trivial, but not Kawamata log terminal.
\end{Def}

Suppose that a log pair $(X,B)$ is log canonical. Then, using log
canonical thresholds, we can compare special tigers for $K_X+B$.
To this end, we introduce
\begin{Def}\com
Let $(X,B)$ be a log canonical pair with nonempty $|-(K_X+B)|$.
The total log canonical threshold of $(X,B)$ is the real number
\[totallct(X,B)=\operatorname{sup}\{r : K_X+B+rD
\textnormal{ is log canonical for any }D\in|-(K_X+B)|\}.\]
\end{Def}
Note that $0\leq totallct(X,B)\leq 1$. If $B=0$, the total log
canonical threshold of $(X,B)$ will be denoted by $totallct(X)$
instead of $totallct(X,B)$. The total log canonical threshold of
$(X,B)$ measures how bad elements of $|-(K_X+B)|$ can be. It is
worthwhile to pay attention to  special tigers realizing the
total log canonical threshold.

\begin{Def}\com
Let $X$ be a normal variety. Let $B$ be an effective
$\mathbb{Q}$-divisor on $X$. A divisor $D$ on $X$ is called a
wild tiger of $K_X+B$ if the following conditions are satisfied:
\begin{enumerate}
\item
$K_X+B+D$ is linearly trivial.
\item
$lct(D;X,B)=totallct(X,B)$.
\end{enumerate}
By a wild tiger for $X$, we mean a wild tiger for $K_X$.
\end{Def}

Comparing with the definition of a special tiger, we note that
numerical triviality is replaced by linear triviality in the
definition of a wild tiger. The concept of a wild tiger is a sort of
counterpart of 1-complement of $K_X+B$. We can find beautiful
interactions between 1-complements and wild tigers in \cite{P99}.
In the present paper, we are mainly interested in finding wild
tigers for smooth hypersurfaces of degree $n\geq 3$ in
$\mathbb{P}^n$.
\renewcommand{\thesection}{\large{\arabic{section}.}}
\section{\hspace{-3mm}\large{Eckardt points}}
\renewcommand{\thesection}{\arabic{section}}
An Eckardt point is a point on a smooth cubic surface $\Sigma$ at
which three lines on $\Sigma$ intersect each other.
In other words, it is a point $p$ on $\Sigma$ such that there is
an element in $|-K_{\Sigma}|$ which is a cone with vertex $p$ and
base consisting of three different points.

Now, we investigate smooth del Pezzo surfaces (in particular
cubic surfaces) to find their wild tigers.
During this investigation, we will observe the special feature of
Eckardt points.
When hunting wild tigers, the first step is to calculate total log
canonical thresholds. But, due to the concrete geometric description
of smooth del Pezzo surfaces,  we may find wild tigers of smooth
del Pezzo surfaces by investigating cubic curves and points in
general position on $\mathbb{P}^2$.

\begin{Prop}\com
Let $S$ be a smooth del Pezzo surface of degree $d$. Then we have
the following table:
\begin{center}
\begin{tabular}{|c|c|c|c|}\hline
$d$ & Remarks & $totallct(S)$& Wild tiger\\ \hline \hline
9 &
\begin{minipage}[m]{.4\linewidth}
\vspace{1mm} Fano index 3 \vspace{1mm}
\end{minipage}
& $\frac{1}{3}$ &
\begin{minipage}[m]{.3\linewidth}
\vspace{1mm} Triple line \vspace{1mm}
\end{minipage}
\\ \hline 8 &
\begin{minipage}[m]{.4\linewidth}
\vspace{1mm} Fano index 1 \vspace{1mm}
\end{minipage}
& $\frac{1}{3}$ &
\begin{minipage}[m]{.3\linewidth}
\vspace{5mm}
\begin{center}
\begin{picture}(50,5)
\multiput(7,0)(30,0){2}{\circle*{4}}
\put(7,0){\line(1,0){30}}
\put(0,5){\textnormal{{\small 3(0)}}}
\put(27,5){\textnormal{{\small 2(-1)}}}
\end{picture}
\end{center}
\vspace{1mm}
\end{minipage}      \\
\hline 8 &
\begin{minipage}[m]{.4\linewidth}
\vspace{1mm}
Fano index 2 \vspace{1mm}
\end{minipage}
 & $\frac{1}{2}$ &
 \begin{minipage}[m]{.3\linewidth}
\vspace{4mm}
\begin{center}
\begin{picture}(50,5)
\multiput(7,0)(30,0){2}{\circle*{4}}
\put(7,0){\line(1,0){30}}
\put(0,5){\textnormal{{\small 2(0)}}}
\put(27,5){\textnormal{{\small 2(0)}}}
\end{picture}

\vspace{1mm}
\begin{picture}(80,5)
\multiput(7,0)(30,0){3}{\circle*{4}}
\put(7,0){\line(1,0){60}}
\put(0,5){\textnormal{{\small 1(0)}}}
\put(27,5){\textnormal{{\small 2(0)}}}
\put(57,5){\textnormal{{\small 1(0)}}}
\end{picture}
\end{center}
\vspace{1mm}
\end{minipage}         \\ \hline
7 &
\begin{minipage}[m]{.4\linewidth}
\vspace{1mm} \vspace{1mm}
\end{minipage}
 & $\frac{1}{3}$  &
 \begin{minipage}[m]{.3\linewidth}
\vspace{5mm}

\begin{center}
\begin{picture}(80,5)
\multiput(7,0)(30,0){3}{\circle*{4}}
\put(7,0){\line(1,0){60}}
\put(0,5){\textnormal{{\small 2(-1)}}}
\put(27,5){\textnormal{{\small 3(-1)}}}
\put(57,5){\textnormal{{\small 2(-1)}}}
\end{picture}
\end{center}
\vspace{1mm}
\end{minipage}       \\ \hline
6 &
\begin{minipage}[m]{.4\linewidth}
\vspace{1mm}
\vspace{1mm}
\end{minipage}
& $\frac{1}{2}$    &
\begin{minipage}[m]{.3\linewidth}
\vspace{4mm}
\begin{center}
\begin{picture}(80,5)
\multiput(7,0)(30,0){3}{\circle*{4}} \put(7,0){\line(1,0){60}}
\put(0,5){\textnormal{{\small 1(-1)}}}
\put(27,5){\textnormal{{\small 2(0)}}}
\put(57,5){\textnormal{{\small 1(-1)}}}
\end{picture}

\vspace{1mm}
\begin{picture}(110,5)
\multiput(7,0)(30,0){4}{\circle*{4}}
\put(7,0){\line(1,0){90}}
\put(0,5){\textnormal{{\small 1(-1)}}}
\put(27,5){\textnormal{{\small 2(-1)}}}
\put(57,5){\textnormal{{\small 2(-1)}}}
\put(87,5){\textnormal{{\small 1(-1)}}}
\end{picture}

\vspace{8mm}
\begin{picture}(50,5)
\put(22,0){\line(0,1){20}}
\put(22,0){\line(-1,-1){14}}
\put(22,0){\line(1,-1){14}}
\put(22,20){\circle*{4}}
\put(22,0){\circle*{4}}
\put(7,-14){\circle*{4}}
\put(36,-14){\circle*{4}}
\put(0,-24){\textnormal{{\small 1(-1)}}}
\put(30,-24){\textnormal{{\small 1(-1)}}}
\put(15,24){\textnormal{{\small 1(0)}}}
\put(27,-2){\textnormal{{\small 2(-1)}}}
\end{picture}
\end{center}
\vspace{8mm}
\end{minipage}      \\ \hline 5 &
\begin{minipage}[m]{.4\linewidth}
\vspace{1mm}
\vspace{1mm}
\end{minipage}
 & $\frac{1}{2}$    &
 \begin{minipage}[m]{.3\linewidth}
\begin{center}
\vspace{10mm}
\begin{picture}(50,5)
\put(22,0){\line(0,1){20}}
\put(22,0){\line(-1,-1){14}}
\put(22,0){\line(1,-1){14}}
\put(22,20){\circle*{4}}
\put(22,0){\circle*{4}}
\put(7,-14){\circle*{4}}
\put(36,-14){\circle*{4}}
\put(0,-24){\textnormal{{\small 1(-1)}}}
\put(30,-24){\textnormal{{\small 1(-1)}}}
\put(15,24){\textnormal{{\small 1(-1)}}}
\put(27,-2){\textnormal{{\small 2(-1)}}}
\end{picture}
\vspace{10mm}
\end{center}
\end{minipage}      \\ \hline
4 &
\begin{minipage}[m]{.4\linewidth}
\vspace{1mm}
\vspace{1mm}
\end{minipage}
& $\frac{2}{3}$      &
\begin{minipage}[m]{.3\linewidth}
\vspace{10mm}
\begin{center}
\begin{picture}(50,5)
\put(7,0){\line(1,0){30}}
\put(7,0){\line(2,3){15}}
\put(37,0){\line(-2,3){15}}
\multiput(7,0)(30,0){2}{\circle*{4}}
\put(22,21){\circle*{4}}
\put(0,-10){\textnormal{{\small 1(-1)}}}
\put(29,-10){\textnormal{{\small 1(-1)}}}
\put(15,24){\textnormal{{\small 1(0)}}}
\end{picture}
\end{center}
\vspace{3mm}
\end{minipage}    \\ \hline 3 &
\begin{minipage}[m]{.4\linewidth}
\vspace{1mm}
$S$ has an Eckardt point. \vspace{1mm}
\end{minipage}
& $\frac{2}{3}$       &
\begin{minipage}[m]{.3\linewidth}
\vspace{10mm}
\begin{center}
\begin{picture}(50,5)
\put(7,0){\line(1,0){30}}
\put(7,0){\line(2,3){15}}
\put(37,0){\line(-2,3){15}}
\multiput(7,0)(30,0){2}{\circle*{4}}
\put(22,21){\circle*{4}}
\put(0,-10){\textnormal{{\small 1(-1)}}}
\put(29,-10){\textnormal{{\small 1(-1)}}}
\put(15,24){\textnormal{{\small 1(-1)}}}
\end{picture}
\end{center}
\vspace{3mm}
\end{minipage}   \\ \hline 3 &
\begin{minipage}[m]{.4\linewidth}
\vspace{1mm}
Generic case \vspace{1mm}
\end{minipage}
& $\frac{3}{4}$        &
\begin{minipage}[m]{.3\linewidth}
\vspace{1mm} A line and a conic intersecting tangentially with
intersection number 2 \vspace{1mm}
\end{minipage} \\ \hline 2 &
\begin{minipage}[m]{.4\linewidth}
\vspace{1mm}
$S$ has an effective anticanonical\\ divisor with a tacnode.
\vspace{1mm}
\end{minipage}
& $\frac{3}{4}$  &
\begin{minipage}[m]{.3\linewidth}
\vspace{1mm} Two lines intersecting tangentially with
intersection number 2\vspace{1mm}
\end{minipage}\\ \hline 2 &
\begin{minipage}[m]{.4\linewidth}
\vspace{1mm}
Generic case \vspace{1mm}
\end{minipage}
& $\frac{5}{6}$   &
\begin{minipage}[m]{.3\linewidth}
\vspace{1mm}
Cuspidal cubic\vspace{1mm}
\end{minipage} \\ \hline 1 &
\begin{minipage}[m]{.4\linewidth}
\vspace{1mm}
\vspace{1mm}
\end{minipage}
& $\frac{5}{6}$    &
\begin{minipage}[m]{.3\linewidth}
\vspace{1mm} Cuspidal cubic\vspace{1mm}
\end{minipage}      \\ \hline
\end{tabular}
\end{center}
where $\bullet$'s denote smooth rational curves, numbers are
multiplicities, and the numbers in parentheses are
self-intersection numbers.
\end{Prop}
\pf The proof is straightforward. Refer to \cite{P99} for some
direction.\qed

We now pay attention to smooth cubic surfaces. We observe that in
this case, there are two different wild tigers. If a cubic
surface has an Eckardt point, then its wild tiger consists of
three lines intersecting at a single point which is an Eckardt
point. If not, then its wild tiger consists of a line and a conic
intersecting tangentially with intersection number 2. Note that
having an Eckardt point is a codimension one condition.

It is expected that Eckardt points indicate special features to del
Pezzo fibrations of degree 3. The evidence can be found in
\cite{Co96} and \cite{P99}. For this reason, it is necessary to
generalize Eckardt points to the case of smooth hypersurfaces of degree
$n\geq 4$ in $\mathbb{P}^n$.

\begin{Def}\label{def}\com
Let $X$ be a smooth hypersurface of degree $n\geq 3$ in $\mathbb{P}^n$.
A point $p$ on $X$ is called an (generalized)
Eckardt point if there is an element $S$ in $|-K_X|$ which is a
cone in $\mathbb{P}^{n-1}$ over a smooth hypersurface of degree $n$
in $\mathbb{P}^{n-2}$ with vertex $p$.
\end{Def}
It is clear that a generalized Eckardt point coincides with the
classical one when $n=3$.

As we noted earlier, if a smooth cubic surface has an Eckardt point,
then its wild tiger is a cone over three different points. This fact
can be generalized as follows:
\begin{Thm}\label{Eckardt}\com
Let $X$ be a smooth hypersurface of degree $n\geq 3$ in
$\mathbb{P}^n$. If $X$ has an Eckardt point $p$, then $S\in
|-K_X|$ as used in Definition~\ref{def} is a wild tiger for $X$.
In particular, $lct(S;X)=totallct(X)=\frac{n-1}{n}$ and the locus
of log canonical singularities $LCS(X,\frac{n-1}{n}S)$ of
$(X,\frac{n-1}{n}S)$ is one point $\{p\}$.
\end{Thm}
We will prove the theorem in the next section.

We see that if a smooth cubic surface has total log canonical
threshold $\frac{2}{3}$, then it has an Eckardt point. It is
natural to conjecture
\begin{Guess}\label{conjecture}\com
Let $X$ be a smooth hypersurface of degree $n\geq 3$ in $\mathbb{P}^n$.
If the total log canonical threshold of $X$ is
$\frac{n-1}{n}$, then a wild tiger $S$ for $X$ is a cone in
$\mathbb{P}^{n-1}$ over a smooth hypersurface of degree $n$ in
$\mathbb{P}^{n-2}$ with vertex $p$. Moreover,
$LCS(X,\frac{n-1}{n}S)=\{p\}$, where $p$ is an Eckardt point of
$X$.
\end{Guess}

Of course, the conjecture holds for $n=3$. In the fourth section,
we will prove the conjecture under the assumption of the Log
minimal model program. Since the Log minimal model program holds
up to dimension $3$, the conjecture is true up to $n=4$.

Before proceeding, we should mention one more problem related to
total log canonical thresholds and wild tigers of smooth Fano
varieties. We see that a generic cubic surface has total log
canonical threshold $\frac{3}{4}$ and that its wild tiger consists
of a line and a conic intersecting tangentially with intersection
number 2. A generic del Pezzo surface of degree 2 has total log
canonical threshold $\frac{5}{6}$, and its tiger is a cuspidal
cubic. This seems to be quite a common phenomenon. So, we ask
what the total log canonical threshold of a generic hypersurface
of degree $n$ in $\mathbb{P}^n$ is and what its wild tiger is.

\renewcommand{\thesection}{\large{\arabic{section}.}}
\section{\hspace{-3mm}\large{A Lower Bound for Total Log Canonical Thresholds}}
\renewcommand{\thesection}{\arabic{section}}
 Let $W$ be a smooth hypersurface of degree $m$ in
$\mathbb{P}^n$ and $H$ be a hyperplane section of $W$, where
$n\geq 4$. It follows from the Lefschetz theorem that the Picard group
of $W$ is a free abelian group generated by a hyperplane section
$H$, i.e., $\operatorname{Pic}(W)=\mathbb{Z}H$. Therefore, a
hyperplane section $H$ of $W$ is irreducible and reduced.

\begin{Lem}\com\label{Pu}
For any curve $C$ on $H$, $\mult_{C}H=1$.
\end{Lem}
\pf Let $p$ be a general point in $\mathbb{P}^{n}\backslash W$.
We consider a cone $P_p$ with the vertex $p$ and the base $C$.
Then we have
\[P_p\cap W=C\cup R_p,\]
where $R_p$ is the residual curve of degree $(m-1)\deg (C)$.
The curves $C$ and $R_p$ intersect at $(m-1)\deg (C)$ different
points (see \cite{Pu95}). Since $H$ is a hyperplane section of
$W\subset \mathbb{P}^n$, we have
\[H\cdot R_p=\deg(R_p)=(m-1)\deg (C).\]
On the other hand,
\[H\cdot R_p\geq \deg(R_p)\mult_CH=(m-1)\deg(C)\mult_CH.\]
\qed

\begin{Cor}\com
A hyperplane section $H$ has only isolated singularities. In
particular, it is normal.
\end{Cor}
\pf The first statement immediately follows from Lemma~\ref{Pu}.
Since $H$ is a smooth in codimension 1 hypersurface of a smooth
variety, it is normal. \qed

\begin{Thm}\label{cheltsov}\com
The log canonical threshold of $H$ in $W$ is
at least $\lambda=\min\{\frac{n-1}{m},1\}$.
\end{Thm}
\pf Let $0<\alpha <\lambda$. We may consider the log pair
$(\mathbb{P}^{n-1},\alpha H) $ instead of $(W,\alpha H)$. Suppose
that $K_{\mathbb{P}^{n-1}}+\alpha H$ is not Kawamata log terminal.
Then the log canonical singularity subscheme
$\mathcal{L}=\mathcal{LCS}(\mathbb{P}^{n-1},\alpha H)$ associated to
$(\mathbb{P}^{n-1}, \alpha H)$ is a
zero-dimensional subscheme. Now, we consider a Cartier divisor
$D$ which is numerically equivalent to
$K_{\mathbb{P}^{n-1}}+\alpha H+(\lambda-\alpha)H'$, where $H'$
is a generic element in $|H|$. Note that
\[\mathcal{O}_{\mathbb{P}^{n-1}}(D)=\left\{\begin{array}{ll}
\mathcal{O}_{\mathbb{P}^{n-1}}(-1) & \mbox{if $m-n\geq -1$}\\
\mathcal{O}_{\mathbb{P}^{n-1}}(m-n)& \mbox{otherwise}
\end{array}
\right. \] By Shokurov's vanishing theorem (see \cite{Am99}),
we have an exact
sequence
\[H^0(\mathbb{P}^{n-1},\mathcal{O}_{\mathbb{P}^{n-1}}(D))\longrightarrow
H^0(\mathcal{L},\mathcal{O}_{\mathcal{L}}(D))\longrightarrow 0.\]
But the first term is zero even though the second term is not.
This is a contradiction. \qed

\begin{Cor}\label{ivan}\com
Suppose that $n=m$. Then the total log canonical  threshold of
$W$ is at least $\frac{n-1}{n}$.
\end{Cor}
\pf This immediately follows from Theorem~\ref{cheltsov}. \qed
\\
{\bfseries\itshape Proof of Theorem~\ref{Eckardt}.} It is obvious
that the log canonical threshold of $S$ with respect to $K_X$ is
$\frac{n-1}{n}$. The theorem then follows from
Corollary~\ref{ivan}.\qed

\begin{Prop}\com\label{one}
If the log canonical threshold $\alpha$ of $H$ in $W$ is not 1,
then the locus of log canonical singularities $LCS(W,\alpha H)$ of
$(W,\alpha H)$ consists of a single point.
\end{Prop}
\pf The proof is similar to that of Theorem~\ref{cheltsov}. \qed

\begin{Ex}\com
\textnormal{Suppose that a hypersurface $W$ in $\mathbb{P}^n$ is given by
equation
\[x_0^m-x_1^m+\sum_{i=2}^{n}x_i^m=0\]
and a hyperplane section $H$ is given by $x_0-x_1=0$. Then the log
canonical threshold of $H$ is $\lambda$. Thus, our $\lambda$ is
the sharp lower bound for log canonical thresholds of hyperplane
sections of smooth hypersurfaces of degree $m$ in $\mathbb{P}^n$.}
\end{Ex}


\renewcommand{\thesection}{\large{\arabic{section}.}}
\section{\hspace{-3mm}\large{Proof of the Conjecture via the
Log Minimal Model Program}}
\renewcommand{\thesection}{\arabic{section}}
Let $X$ be a smooth hypersurface of degree $n\geq 4$ in $\mathbb{P}^n$.
Let $S$ be a hyperplane section of $X$. The goal
of this section is to prove Conjecture~\ref{conjecture} for $n=4$.
Specifically, under the assumption of the Log minimal model
program in dimension $\leq n-1$, we will show that if a log pair
$(X,\frac{n-1}{n}S)$ is not Kawamata log terminal, then $S$ is a
cone in $\mathbb{P}^{n-1}$ over a smooth hypersurface of degree
$n$ in $\mathbb{P}^{n-2}$.

We suppose that $(X, \frac{n-1}{n}S)$ is not Kawamata log
terminal. Then, $LCS(X, \frac{n-1}{n}S)$ is nonempty and consists
of only finite number of points (in fact, only one point by
Proposition~\ref{one}). So, we may forget about the hypersurface $X$
and deal only with the log pair $(\mathbb{P}^{n-1},
\frac{n-1}{n}S)$.

From now on, we assume that the Log minimal model program  holds
for dimension $\leq n-1$.
\begin{Lem}\com
There exists a birational morphism $f:V\longrightarrow
\mathbb{P}^{n-1}$ satisfying the following:
\begin{itemize}
\item
$f$ is an isomorphism outside of $LCS(\mathbb{P}^{n-1},
\frac{n-1}{n}S)=\{p\}$,
\item
$V$ has $\mathbb{Q}$-factorial terminal singularities, and
\item
there is an effective $f$-exceptional $\mathbb{Q}$-divisor $E$ on
$V$ such that the support of $E$ coincides with that of the
$f$-exceptional locus, $\lfloor E\rfloor \ne 0$, and
$K_V+\frac{n-1}{n}f^{-1}_*S+E=f^*(K_{\mathbb{P}^{n-1}}+\frac{n-1}{n}S)$.
\end{itemize}
\end{Lem}
\pf Let $g:V'\longrightarrow \mathbb{P}^{n-1}$ be a log terminal
blow-up of $(\mathbb{P}^{n-1},\frac{n-1}{n}S)$. Since $V'$ has
$\mathbb{Q}$-factorial Kawamata log terminal singularities, we
may take a terminal blow-up $h:V\longrightarrow V'$ with respect
to $V'$ (see \cite{Pro99b}). Then the birational morphism
$f=g\circ h:V\longrightarrow \mathbb{P}^{n-1}$ will satisfy the
conditions. \qed

We fix such a birational morphism
$f:V\longrightarrow\mathbb{P}^{n-1}$. Let $\tilde{S}=f^{-1}_*(S)$.
Since the log pair $(\mathbb{P}^{n-1}, \frac{n-1}{n}S)$ is log
canonical, the log pair $(V,\frac{n-1}{n}\tilde{S}+E)$ is also log
canonical. We see that $K_V+\frac{n-1}{n}\tilde{S}+E$ is not nef
and that $-(K_V+\frac{n-1}{n}\tilde{S}+E)$ is not ample.
Therefore, there is an extremal contraction $g:V\longrightarrow
W$ such that $-(K_V+\frac{n-1}{n}\tilde{S}+E)$ is $g$-ample.
Because $-(K_V+\frac{n-1}{n}\tilde{S}+E)$ is $g$-ample and
$f$-numerically trivial, no curve contracted by $g$ is contained
in the fibers of $f$. In particular, $W$ is not a point.

\begin{Lem}\label{ll}\com
Suppose that the extremal contraction $g$ contracts a subvariety $F$
of $V$ to a subvariety $Z$ of $W$. Then $\dim F-\dim Z=1$.
\end{Lem}
\pf Suppose that $\dim F-\dim Z> 1$. Let $G$ be a fiber of $g$
over $Z$. Since $G\cap E\ne \emptyset$ and $V$ is
$\mathbb{Q}$-factorial, there is a curve on $G\cap E$ which is
contracted by both $f$ and $g$. But this is impossible. \qed

\begin{Prop}\com
If the extremal contraction $g:V\longrightarrow W$ is a Mori fiber
space, then $g$ is a conic bundle.
\end{Prop}
\pf This  follows from Lemma~\ref{ll}.
\qed\\
The following lemma is due to V.~Shokurov's paper (\cite{Sho01}),
which is in
preparation. It generalizes X. Benveniste and S.
Mori's results (\cite{Ben85} and \cite{Mo88}) under the assumption
of the existence of flips.

\begin{Lem}\label{vshokurov}\com
Suppose that $Y$ has at worst $\mathbb{Q}$-factorial terminal
singularities. Let $h:Y\longrightarrow Z$ be a birational
contraction. If a curve $C$ on $Y$ is an irreducible component of
the exceptional locus of $h$, then $K_Y\cdot C>-1$.
\end{Lem}
\pf It is enough to consider the statement over an analytic
neighborhood of $h(C)=q\in Z$. We choose a divisor $H$ on $Y$ with
$H\cdot C=1$. In addition, we may assume that the exceptional
locus of $h$ is the curve $C$. Suppose that $K_Y\cdot C\leq -1$.
Then $(K_Y+H)\cdot C=0$ (see \cite{Sho01}). We consider an $K_Y$-flip of $h$;
\[ \xymatrix{
Y \ar@{-->}[rr]^{\phi}  \ar[dr]_h && Y^+ \ar[dl]^{h^+}\\
&Z}\]
Let $H^+$ be the birational transform of $H$ to $Y^+$ via $\phi$.
Since we have $\dim Ex(h) +\dim Ex(h^+)\geq \dim Y-1$ (see
\cite{KMM}), we get $\dim Ex(h^+)= \dim Y-2$.

Note that the numerical $h$-triviality of $K_Y+H$  implies $H^+\cdot
C'<0$ for any curve $C'$ on $Ex(h^+)$. Therefore, $Ex(h^+)\subset
H^+$.

Let $E$ be the exceptional divisor of blow-up centered at a
component of $Ex(h^+)$. We may assume that the center of $E$ on
$Y$ is not contained in $H$. Then we have the following
inequality:
\[a(E;Y,H)\leq a(E;Y^+,H^+)\leq 0,\]
where $a(E;Y,H)$ (resp. $a(E;Y^+,H^+)$) is the discrepancy of $E$
with respect to $K_Y+H$ (resp. $K_{Y^+}+H^+$). This is a
contradiction because $Y$ is terminal.\qed

\begin{Lem}\com
The extremal contraction $g$ is not a small contraction whose
exceptional locus has a curve as an irreducible component.
\end{Lem}
\pf Suppose that $g$ is a small contraction whose exceptional locus
has a curve $C$ as an irreducible component. We have $K_V\cdot
C>-1$ by Lemma~\ref{vshokurov}. On the other hand,
\[(K_V+\frac{n-1}{n}\tilde{S})\cdot C
=(f^*\mathcal{O}_{\mathbb{P}^{n-1}}(-1)-E)\cdot C=-\deg(f_*C)-E\cdot C
\leq -1.\] Thus,
$\tilde{S}\cdot C<0$ and $C\subset \tilde{S}$, and hence
\[(K_V+\tilde{S})\cdot C<(K_V+\frac{n-1}{n}\tilde{S})\cdot C\leq -1.\]

Let $\nu:\hat{S}\longrightarrow\tilde{S}$ be a normalization of
$\tilde{S}$. By adjunction, we have
\[K_{\hat{S}}+\Diff_{\tilde{S}}(0)=\nu^*((K_V+\tilde{S})|_{\tilde{S}}).\]
Since $\tilde{S}$ is smooth at a generic point of the curve $C$,
the curve $\nu^{-1}_*C$ cannot be contained in
$\Diff_{\tilde{S}}(0)$, and hence
$K_{\hat{S}}\cdot\nu^{-1}_*C<-1$. On the other hand,
$K_{\hat{S}}\cdot\nu^{-1}_*C\geq-1$ since the curve $\nu^{-1}_*C$
is contractible. \qed

So far, we proved that the extremal contraction $g$ is a contraction
of a subvariety $F$ of $V$ to a subvariety $Z$ of $W$ with $\dim
F-\dim Z=1$ and $\dim Z\geq 1$. Let $C$ be a general fiber of
the morphism $g$ over $Z$. Then we have
\[(K_V+\frac{n-1}{n}\tilde{S})\cdot C=
(f^*\mathcal{O}_{\mathbb{P}^{n-1}}(-1)-E)\cdot C =
-\deg(f_*C)-E\cdot C <-\deg(f_*C),\] because the curve $C$ should
meet the exceptional locus of $f$.

\begin{Lem}\label{done}\com
If the extremal contraction $g:V\longrightarrow W$ is not a conic
bundle, then $F=\tilde{S}$.
\end{Lem}
\pf
Let $C$ be a general enough fiber of the morphism $g$.
The inequality
\[-1\leq
K_V\cdot C= (K_V+\frac{n-1}{n}\tilde{S})\cdot
C-\frac{n-1}{n}\tilde{S}\cdot C
<-\deg(f_*C)-\frac{n-1}{n}\tilde{S}\cdot C\] implies that
$F\subset\tilde{S}$. If the codimension of the subvariety $F$ of $V$
is greater than 1, then we get
\[-1\geq -\deg(f_*C)>(K_V+\frac{n-1}{n}\tilde{S})\cdot C =
(K_V+\tilde{S})\cdot C-\frac{1}{n}\tilde{S}\cdot C>
-1-\frac{1}{n}\tilde{S}\cdot C>-1,\] where the second to the last
inequality is implied by Lemma~\ref{vshokurov}. This is absurd.
Hence, the subvariety $F$ has codimension 1. Consequently, $F=\tilde{S}$.
\qed
\\
Now we know that the extremal contraction $g$ is either a conic
bundle or a contraction of the divisor $\tilde{S}$ of $V$ to a subvariety $Z$
of $W$ with $\dim Z=n-3$.

\begin{Thm}\com\label{cone}
If the Log minimal model program holds for dimension $\leq n-1$,
then Conjecture~\ref{conjecture} holds for $n$.
\end{Thm}
\pf Suppose that $g$ is a conic bundle. Then $\lfloor
E\rfloor\cap C \ne \emptyset $ since no component of the divisor $E$
lies in the fibers of $g$. Therefore,
\[\deg(f_*C)=-(K_V+\frac{n-1}{n}\tilde{S}+E)\cdot
C=2-\frac{n-1}{n}\tilde{S}\cdot C-E\cdot C<2.\] Consequently,
$f_*C$ is a line on $S$. This implies that $S$ is a cone in
$\mathbb{P}^{n-1}$.

Now, we suppose that the morphism $g$ is not a conic bundle. Then it
is a contraction of the divisor $\tilde{S}$ of $V$ to a subvariety $Z$ of $W$
with $\dim Z=n-3$.
Therefore, we have
\[-\deg(f_*C)>(K_V+\frac{n-1}{n}\tilde{S})\cdot C =
(K_V+\tilde{S})\cdot C-\frac{1}{n}\tilde{S}\cdot C=
-2-\frac{1}{n}\tilde{S}\cdot C>-2.\] Thus, $f_*C$ is a line on
$S$. Consequently, $S$ is a cone in $\mathbb{P}^{n-1}$. \qed
\begin{Cor}\label{cor}\com
Conjecture~\ref{conjecture} holds for $n=4$.
\end{Cor}
\pf The Log minimal model program has been proven for dimension 3.
Theorem~\ref{cone} implies the statement. \qed
\begin{Cor}\com
The total log canonical  threshold of a smooth quartic $X$ in
$\mathbb{P}^4$ is $\frac{3}{4}$ if and only if the quartic $X$
has an Eckardt point.
\end{Cor}
\pf This immediately follows from Corollary~\ref{cor}. \qed

\renewcommand{\thesection}{\large{\arabic{section}.}}
\section{\hspace{-3mm}\large{
Application}}
\renewcommand{\thesection}{\arabic{section}}
Let $\mathcal{O}$ be a discrete valuation ring with quotient field $K$.
We assume that the residue
field is of characteristic zero.
For a scheme $\pi : X\longrightarrow \spec \mathcal{O}$,
we denote its scheme-theoretic fiber
$\pi^*(o)$ by $S_X$, where $o$ is the closed point of
$\spec \mathcal{O}$.

\begin{Thm}\label{sunjoo}\com
Let $X$ and $Y$ be smooth Fano fibrations over $\spec \mathcal{O}$
such that
their generic fibers are each isomorphic to a smooth hypersurface of
degree $n$ in $\mathbb{P}^{n}_K$, where $n\geq 3$.
Then any birational map of $X$ into $Y$ over $\spec \mathcal{O}$
which is biregular
on the generic fiber is biregular.
\end{Thm}
\pf Note that our birational map cannot be an isomorphism in codimension 1
(see \cite{Co95}).
Anticanonical divisors of $S_X$ and $S_Y$ are
very ample. Moreover, their total log canonical
thresholds are strictly larger than
$\frac{1}{2}$.
Therefore, the same method as in \cite{P99} works.
\qed
\begin{Rmk}
\textnormal{Birational rigidity has been proven for any smooth
hypersurface of degree $n$ in $\mathbb{P}^n$ and any generic
hypersurface of degree $m$ of $\mathbb{P}^m$, where $4\leq n\leq
8$ and $m\geq 9$. For details, refer to \cite{Che00}, \cite{Is79},
\cite{IsMa71}, \cite{Pu87}, and \cite{Pu98}. Using birational
rigidity, we can easily prove Theorem~\ref{sunjoo} in the case of
$n\leq 8$. Also, we can obtain a weaker statement than
Theorem~\ref{sunjoo} for $n\geq 9$.}
\end{Rmk}
\textbf{Acknowledgments.}  The authors would like to thank Prof.~V.~Shokurov
for many invaluable conversations.

\footnotesize
\nocite{Mu00}
\bibliographystyle{amsplain}
\bibliography{Eck}
\vspace{5mm}
\begin{tabbing}
  \hspace*{29 em}\=\kill
Ivan Cheltsov\>Jihun Park\\
Division of Pure Mathematics\>Department of Mathematics\\
University of Liverpool\>The Johns Hopkins University\\
Peach Street, Liverpool, L69 7ZL, England, U.K.\>
Baltimore, MD. 21218, U.S.A.   \\
\texttt{cheltsov@yahoo.com}\>\texttt{wlog@bigfoot.com}
\end{tabbing}

\end{document}